\documentclass{yfNote}

\yfnotesetup{
  info = {
    title = {A Perturbative Two-Eigenvalue Faber-Krahn Problem for the Gaussian Short-Time Fourier Transform},
    author = {Yifan Zhang},
    date = {}
  },
  hyperref = {
    color = academic
  },
  bib = {
    style = numeric,
    resource = {references.bib}
  }
}
\author{
  Yifan Zhang
}
\yfnotebibsetup{

}
\renewcommand{\thesection}{\arabic{section}}

\counterwithout{theorem}{chapter}
\counterwithin{theorem}{section}
\renewcommand{\thetheorem}{\thesection.\arabic{theorem}}

\renewcommand{\theHtheorem}{\theHsection.\arabic{theorem}}

\makeatletter
\@for\yf@env:=lemma,proposition,corollary,definition,example,remark,note\do{%
  \expandafter\def\csname the\yf@env\endcsname{\thetheorem}%
  \expandafter\def\csname theH\yf@env\endcsname{\theHtheorem}%
}
\makeatother

\numberwithin{equation}{section}

\newcommand{\norm}[1]{\left\lVert #1 \right\rVert}
\newcommand{\inner}[2]{\left\langle #1, #2 \right\rangle}
\usepackage{mathtools,graphicx}

\setlist[enumerate,1]{label=(\roman*)}
\AtBeginDocument{%
}
\declaretheorem[
  name=Conjecture,
  sibling=theorem,
  style=plain
]{conjecture}
\AtBeginDocument{%
}
\NewCommandCopy{\OriginalSection}{\section}

\RenewDocumentCommand{\section}{s o m}{%
  \IfBooleanTF{#1}{%
    \OriginalSection*{#3}%
    \addcontentsline{toc}{section}{#3}%
  }{%
    \IfNoValueTF{#2}
      {\OriginalSection{#3}}
      {\OriginalSection[#2]{#3}}%
  }%
}
\begin{document}

\maketitle
% 摘要
\begingroup
  \normalfont\fontsize{10pt}{12pt}\selectfont
  \setlength{\leftskip}{3pc}
  \setlength{\rightskip}{3pc}
  \setlength{\parindent}{0pt}
  \textsc{Abstract.}\enspace
  We prove that, for every fixed $s\in(0,\infty)$, the set of area $s$ maximizing the weighted sum of the first two eigenvalues of the associated Toeplitz operators on the Fock space must be a disc, up to translations and sets of measure zero. To do so, we first prove the existence of maximizer sets. Then we show that every maximizer set must be a perturbation of a disc when the weight is sufficiently small. Finally, we use the vector-valued version of the analytic implicit function theorem and the inverse function theorem to conclude that the perturbation must vanish when the weight tends to zero.
  \par
\endgroup

% 首页底部的无编号关键词
\begingroup
  \renewcommand{\thefootnote}{}
  \footnotetext[0]{%
  }
\endgroup
\tableofcontents

\section {Introduction}\label{intro}
\subsection{Main results}
The Fock space, denoted by $\mathcal{F}$, is the Hilbert space of all entire functions $f:\mathbb{C} \to \mathbb{C} $ such that \[\norm{f}^2=\int_\mathbb{C} |f(z)|^2e^{-\pi|z|^2} \dd A(z)<\infty,\]where $\dd A(z)$ denotes the Lebesgue measure on $\mathbb{C} $, with the natural inner product \[\inner{f}{g}=\int_\mathbb{C} f(z)\overline{g(z)}e^{-\pi|z|^2} \dd A(z).\]Let $P: L^2(\mathbb{C}, e^{-\pi|z|^2}\dd A(z))\to\mathcal{F} $ denote the orthogonal projection onto the Fock space, the Toeplitz operator associated with a measurable set $\Omega\subset \mathbb{C} $ is defined as \[T_{\Omega}f(z) :=P(\textbf{1}_{\Omega}f)(z), \quad z\in\mathbb{C} \quad\text{and} \quad f\in\mathcal{F}.\] By standard arguments, $T_{\Omega}$ is a positive operator in $\mathfrak{S}_1$, the Schatten 1-class on $\mathcal{F} $, whenever $0<|\Omega|<\infty$. Here, $|\Omega|$ denotes the Lebesgue measure of $\Omega$ and $\mathfrak{S} _1$ is a Banach space with the norm \[\norm{T}_{\mathfrak{S} _1} := \sum\limits_{n=0}^{\infty}\inner{|T|e_n}{e_n}.\]Indeed, we have 
\begin{equation}\label{s1}
\norm{T_{\Omega}}_{\mathfrak{S} _1}=|\Omega|.
\end{equation}
Therefore, $T_{\Omega}$ is a compact operator on $\mathcal{F} $, which admits a nondecresing sequence of positive eigenvalues denoted by \[1>\lambda_1(\Omega)\geq\lambda_2(\Omega)\geq\cdots>0.\]The main result of \cite{MR4480144} can be stated as follows:
\begin{theorem}\label{nt22}
 For every measurable set $\Omega\subset \mathbb{C} $ satisfying $0<|\Omega|<\infty$, we have
 \begin{equation}\label{eqnt22}
  \lambda_1(\Omega)\leq\lambda_1(\Omega^{\ast}),
 \end{equation}  
where $\Omega^{\ast}$ denotes the disc with the same area as $\Omega$. Moreover, equality is attained if and only if $\Omega$ coincides (up to a set of measure zero) with a disc.
\end{theorem}   
This naturally raises the following conjecture:
\begin{conjecture}\label{con}
For every $N\geq1$ and measurable set $\Omega\subset \mathbb{C} $ satisfying $0<|\Omega|<\infty$, it holds
\begin{equation}\label{eqcon}
  \sum\limits_{n=1}^{N}\lambda_n(\Omega)\leq\sum\limits_{n=1}^{N}\lambda_n(\Omega^{\ast}).
 \end{equation}
 Moreover, equality is attained if and only if $\Omega$ coincides (up to a set of measure zero) with a disc.
\end{conjecture}   
In \cite{MR5109974} it was proved that \eqref{eqcon} holds for all radially symmetric sets. More recently, Nicola, Riccardi, and Tilli \cite{nicola2026existenceextremizerssumeigenvalues} proved the existence of a maximizer set such that the supremum 
\[
S_N(s) :=\sup\limits_{|\Omega|=s}\sum\limits_{n=1}^{N}\lambda_n(\Omega)
\]
is attained. Later, Abreu \cite{abreu2026isoperimetricinequalitykyfan} proved that $S_N(s)$ can be attained by the disc with the same area as $\Omega$. However, the uniqueness (modulo translations and sets of measure zero) of the maximizer set remains open. In \cite{MR5109974}, the authors obtained some negative results on this problem. It was shown that the maximizer set of $\sup\limits_{|\Omega|=s}\lambda_2(\Omega)$ can only be an annulus when $\Omega$ is radially symmetric. This conclusion fails if the radial symmetry assumption is dropped and area $s$ is small enough. 
\vspace{0.3em}  
\par This paper gives an affirmative answer in part to this problem. Our result reads as follows:
\begin{theorem}\label{main}
For every $s>0$, there exists $\delta_s>0$ such that, whenever $0<\delta<\delta_s$, every maximizer of 
\begin{equation}\label{main def}
J_{\delta}(\Omega) := \lambda_1(\Omega)+\delta\lambda_2(\Omega)
\end{equation}
among measurable sets satisfying $|\Omega|=s$ is a disc up to translation.
\end{theorem}  
\begin{remark}\label{coro}
Notice that \autoref{main} can easily be proved once \autoref{con} holds. Indeed, setting $N=2$ and $\delta_s=1$, we obtain that 
\[
\lambda_1(\Omega)+\lambda_2(\Omega)\leq \lambda_1(\Omega^{\ast})+\lambda_2(\Omega^{\ast})
\]
and hence 
\[\begin{aligned}
\lambda_1(\Omega)+\delta\lambda_2(\Omega)&=(1-\delta)\lambda_1(\Omega)+\delta(\lambda_1(\Omega)+\lambda_2(\Omega))\\ &\leq(1-\delta)\lambda_1(\Omega^{\ast})+\delta(\lambda_1(\Omega^{\ast})+\lambda_2(\Omega)^{\ast})\\ &=\lambda_1(\Omega^{\ast})+\delta\lambda_2(\Omega^{\ast})
\end{aligned}\]
by \eqref{eqnt22}. The equality holds if and only if 
\[
\lambda_1(\Omega)= \lambda_1(\Omega^{\ast})
\]
and 
\[
\lambda_1(\Omega)+\lambda_2(\Omega)=\lambda_1(\Omega^{\ast})+\lambda_2(\Omega)^{\ast},
\]
which yields the desired conclusion.
\end{remark}  
\subsection{Outline and notation of the paper}
The paper is organized as follows. In \autoref{pre}, we give some notations and preliminary results on Fock spaces. In \autoref{exi} we prove the existence of maximizer sets attaining the supremum in \eqref{main def}. In \autoref{geom} we give a geometric description of super-level sets when $\delta$ is sufficiently small. Finally, in \autoref{analy} we prove \autoref{main}.
\vspace{0.3em}
\par We conclude with some notation used throughout the paper. For $A,B\in\mathbb{R}$, the symbol $A\lesssim B$ means that there exists $C>0$ such that $A\leq CB$. If we allow $C$ to explicitly depend on a parameter $p$, we will write $A\lesssim_pB$. Additionally, we will use $\mathbb{T}$ to denote the unit circle, which we identify
with the compact Abelian group $\mathbb{R}/2\pi\mathbb{Z}=[-1/2,1/2)$. Given an integrable function $f:\mathbb{T} \to \mathbb{C}$, we define its Fourier transform by 
\[
\widehat{f}(n) :=\int_{\mathbb{T} }f(\theta)e^{-in\theta}\ \dd\theta, \quad n\in\mathbb{Z}.
\]
\section{Preliminary results}\label{pre}
\subsection{Properties of the Fock space}
It is well known that $\mathcal{F} $ is a reproducing Hilbert space. For $w\in\mathbb{C} $ we denote the reproducing kernels by $K_w$ and normalized forms by $k_w$. For $f\in\mathcal{F}$, the Husimi fuction is defined as\[u_f(z) :=|f(z)|^2e^{-\pi|z|^2}, \quad z\in\mathbb{C}.\] More generally, for $\delta>0$ and $F=(f_1,f_2)\in\mathcal{F} ^2$, the $\delta$-joint Husimi function is given by \[u_F(z) :=|f_1(z)|^2e^{-\pi|z|^2}+\delta|f_2(z)|^2e^{-\pi|z|^2}, \quad z\in\mathbb{C}.\]Here, the product $\mathcal{F}^2$ is equipped with $\delta$-weighted inner product \[\inner{F}{G}_{\delta} :=\inner{f_1}{g_1}+\delta\inner{f_2}{g_2},\quad F,G\in\mathcal{F} ^2.\]Finally, the Weyl unitary operators $U_w:\mathcal{F} \to \mathcal{F} $ are defined for $w\in\mathbb{C} $ by \[U_wf(z)=f(z)e^{\pi z\bar{w}}e^{\pi|w|^2/2},\quad z\in\mathbb{C}.\]It is very easy to check that \[u_{U_wf}(z)=u_f(z-w)\]for every $z,w\in\mathbb{C}$. Therefore, $U_w$ are the unitary operators on $\mathcal{F} $. The same is true for the $\delta$-joint Husimi function of a tuple $F$, namely \[u_{U_wF}(z)=u_F(z-w)\]where $F=(f_1,f_2)$ and $U_wF=(U_wf_1,U_wf_2)$. The followings are some basic results on Fock spaces, which we will not prove in this paper. For necessary proofs and more details on Fock spaces, the reader is referred to \cite{MR2934601}.

\begin{lemma}\label{fock lem}
For $f\in\mathcal{F} $, we have
\vspace{1em} 
\begin{enumerate}
\item $K_w(z)=e^{\pi\bar{w}z}, \quad w,z\in\mathbb{C}$.  
\vspace{1em} 
\item $\norm{u_f}_{1}=\norm{f}$.
\vspace{1em} 
\item For every $z\in\mathbb{C} $, \[u_f(z)=|\inner{f}{k_z}|^2\leq \norm{f}^2.\]
\item $\lim\limits_{|z|\to\infty}u_f(z)=0.$
\vspace{1em} 
\item The functions \[e_n(z) :=\Bigl(\frac{{\pi}^n}{n!}\Bigr)^{1/2}z^n\]form an orthonormal basis of $\mathcal{F} $.
\vspace{1em} 
\item For every $z\in\mathbb{C} $, $f$ admits the expansion \[f(z)=\sum\limits_{n=0}^{\infty}a_ne_n(z), \quad \sum\limits_{n=0}^{\infty}|a_n|^2=\norm{f}^2.\]
\end{enumerate}  
\end{lemma}  
\subsection{The weighted Ky Fan reduction}
We now state the following lemma, which is a straightforward application of the Ky Fan Principle. 

\begin{lemma}\label{kf}
For every $\delta\in[0,1]$, it holds \[J_{\delta}(\Omega)=\sup\{\inner{T_{\Omega}F}{F}_{\delta}: F \ \text{is an orthonormal system in} \ {\mathcal{F} }\}. \]
\end{lemma}  
\begin{proof}
On the one hand, by the Ky Fan Principle, for every orthonormal system $F=(f_1,f_2)$ in ${\mathcal{F}} $, we have 
\[\begin{aligned}
\inner{T_{\Omega}F}{F}_{\delta}&=(1-\delta)\inner{T_{\Omega}f_1}{f_1}+\delta(\inner{T_{\Omega}f_1}{f_1}+\inner{T_{\Omega}f_2}{f_2})\\&\leq(1-\delta)\lambda_1(\Omega)+\delta(\lambda_1(\Omega)+\lambda_2(\Omega))\\&=J_{\delta}(\Omega).
\end{aligned}\]
On the other hand, by the standard operator theory, there exists an orthonormal system $\widetilde{F}=(\tilde{f_1},\tilde{f_2})$ such that \[\inner{T_{\Omega}\tilde{f_1}}{\tilde{f_1}}=\lambda_1(\Omega), \quad\inner{T_{\Omega}\tilde{f_1}}{\tilde{f_1}}+\inner{T_{\Omega}\tilde{f_2}}{\tilde{f_2}}=\lambda_1(\Omega)+\lambda_2(\Omega)\]and hence \[J_{\delta}(\Omega)=\inner{T_{\Omega}\widetilde{F}}{\widetilde{F}}_{\delta},\]which yields the desired conclusion.
\end{proof}  
As a corollary of \autoref{kf} and the bathtub principle, we obtain the following:

\begin{corollary}\label{reduction}
For every $s>0$ we have
\[\begin{aligned}
\mathcal{S}_{\delta}(s) :=\sup\limits_{|\Omega|=s}J_{\delta}(\Omega)&=\sup\{\int_{\Omega}u_F(z) \ \dd A(z): F\in{\mathcal{F} }^2 \text{orthonormal system and} \ |\Omega|=s\}\\&=\sup\{\int_{\Omega_F(s)}u_F(z) \ \dd A(z):F\in{\mathcal{F} }^2 \text{orthonormal system}\}.
\end{aligned}\]  
Here $\Omega_F(s)$ denotes the unique super--level set of $u_F$ of measure s, i.e., \[\Omega_F(s)=\{u_F>t\} :=\{z\in\mathbb{C}: u_F(z)>t\},\]where $t>0$ is determined by $|\Omega_F(s)|=s$. Futhermore, if $\mathcal{S}_{\delta}(s)$ could be attained, then the maximizer set $\Omega$ must concide with $\Omega_F(s)$, up to a set of measure 0, for some $F\in{\mathcal{F} }^2$.
\end{corollary}  
\begin{proof}
The first equation follows directly from \autoref{kf}. For the second, it suffices to observe that 
\[\begin{aligned}
\int_{\Omega}u_F&=\int_{\Omega\setminus\Omega_F(s)}u_F+ \int_{\Omega\cap  \Omega_F(s)} u_F\\&\leq t|\Omega\setminus\Omega_F(s)|+ \int_{\Omega\cap  \Omega_F(s)} u_F\\&=t|\Omega_F(s)\setminus\Omega|+ \int_{\Omega\cap  \Omega_F(s)} u_F\\&\leq\int_{\Omega_F(s)\setminus\Omega}u_F+ \int_{\Omega\cap  \Omega_F(s)} u_F\\&=\int_{\Omega_F(s)}u_F.
\end{aligned}\] 
The seconde equation above follows since $|\Omega|=|\Omega_F(s)|=s$.
\vspace{0.3em}
\par Finally, if $\mathcal{S}_{\delta}(s)$ could be attained, then there exist some $\Omega\subset \mathbb{C} $ with $|\Omega|=s$ and $F\in{\mathcal{F} }^2$ such that \[\mathcal{S}_{\delta}(s)=J_{\delta}(\Omega)=\int_{\Omega}u_F.\]Since \[\int_{\Omega}u_F\leq\int_{\Omega_F(s)}u_F\leq\mathcal{S}_{\delta}(s),\]we have \[\int_{\Omega_F(s)}u_F=\mathcal{S}_{\delta}(s).\]By the arguments above, the equality holds if and only if \[|\Omega\setminus\Omega_F(s)|=|\Omega_F(s)\setminus\Omega|=0,\]which completes the proof.
\end{proof}  
\begin{remark}
According to \autoref{fock lem}.4, if a maximizer set exists then, up to a set of measure 0, it must be bounded.
\end{remark}  
\section{The existence of  maximizer sets}\label{exi}
In this section we prove that for every $\delta\in[0,1]$ and $s\in(0,\infty)$, $\mathcal{S} _{\delta}(s)$ is attained. Although we have already proved the existence in \autoref{coro}, we nevertheless give an alternative proof. The methods that we will use in what follows are essentially the same as those in \cite{nicola2026existenceextremizerssumeigenvalues}. 

\subsection{Preliminary lemmas}
In this subsection we gather several lemmas that will be useful in the proofs below. Although some of these lemmas are slightly different from their original forms in \cite{nicola2026existenceextremizerssumeigenvalues}, the essential ideas behind the proofs are the same. For this reason, we omit the proofs here.
\begin{lemma}\label{lip}
Let $F=(f_1,f_2)$, $G=(g_1,g_2)\in{\mathcal{F} }^2$. Then \[|{\mathcal{J}  }_{\delta}(F)-{\mathcal{J}  }_{\delta}(G)|\lesssim_{\delta}(\norm{F}+\norm{G})\norm{F-G}.\] Here ${\mathcal{J}  }_{\delta}(F)$ is defined as \[{\mathcal{J}  }_{\delta}(F) :=\int_{{\Omega}_F(s)}u_F(z) \ \dd A(z).\]
\end{lemma}  
\begin{lemma}\label{unitary}
For every $F=(f_1,f_2)\in{\mathcal{F} }^2$ there exists a unitary matrix $U\in\mathbb{C} ^{2\times2}$ such that, if $G$ is the 2-tuple $(g_1,g_2)\in{\mathcal{F} }^2$ defined by \[
\begin{pmatrix}
  g_1(z)\\[6pt]
  {\delta}^{1/2}g_2(z)
\end{pmatrix}
:= U
\begin{pmatrix}
  f_1(z)\\[6pt]
  {\delta}^{1/2}f_2(z)
\end{pmatrix},
\qquad z\in\mathbb{C},
\]then 
\begin{enumerate}
\vspace{1em}
\item the functions in $G=(g_1,g_2)$ are pairwise orthogonal;
\vspace{1em} 
\item the $\delta$-joint Husimi function is preserved, that is $u_F(z)=u_G(z)$ for every $z\in\mathbb{C} $. 
\end{enumerate}  
\vspace{1em}
Moreover, if $F$ is an orthonormal system in $\mathcal{F} $, then so is $G$ defined aboe, for every unitary matrix $U$.
\end{lemma}
Before starting the next lemma, we introduce a necessary notation. In particular, for every $f\in\mathcal{F} $ we let \[\widehat{f} :=
\begin{cases}
\frac{f}{\norm{f}}, &f\neq0,\\[6pt]
0, &f=0.
\end{cases}\] 
Futhermore, if $F=(f_1,f_2)\in{\mathcal{F} }^2$, we let $\widehat{F}=(\widehat{f_1},\widehat{f_2}).$
\begin{lemma}\label{decom}
Let $F=(f_1,f_2)$, $H=(h_1,h_2)\in{\mathcal{F} }^2$ such that \[\norm{f_k}^2+\norm{h_k}^2=1\]for every $k=1,2$. Assume that $F$ is ordered so that $\norm{f_1}\geq\norm{f_2}$ and let \[a_0=1-\norm{f_1}^2,\quad a_1=\norm{f_1}^2-\norm{f_2}^2,\quad a_2=\norm{f_2}^2.\]We introduce the mixed blocks 
\[\widehat{M}_k=
\begin{cases}
\widehat{H},&k=0,\\  
(\widehat{f}_1,\widehat{h}_2),&k=1\\
\widehat{F},&k=2.
\end{cases}\]  
Then, it holds \[u_F+u_H=\sum\limits_{k=0}^2a_ku_{\widehat{M}_k}.\] 
\end{lemma}  
\begin{lemma}\label{conti}
For every $\delta\in[0,1]$ the function $s\mapsto {\mathcal{S} }_{\delta}(s)$ is continuous.
\end{lemma}  
\begin{lemma}\label{max}
For every $\delta\in(0,1]$ and for every $s\in(0,\infty)$ it holds \[\max\limits_{t\in[0,s]}({\mathcal{S} }_0(t)+\delta{\mathcal{S} }_0(s-t))<{\mathcal{S} }_{\delta}(s).\]
\end{lemma} 
\begin{lemma}\label{hilbert}
Let $\mathcal{H} $ be a Hilbert space with $\dim \mathcal{H} \geq4$, and assume that $x=(x_1,x_2), y=(y_1,y_2)\in{\mathcal{H} }^2$ are such that 
\begin{enumerate}[label=(\alph*)]
\vspace{1em}
\item $x+y$ is an orthonormal system. 
\vspace{1em}
\item $x$ is an orthogonal system with $\norm{x_k}\leq1$ for every $k=1,2$.
\vspace{1em}
\item $|\inner{f_k}{g_l}|<\epsilon$ for some $\epsilon\in(0,1)$ and for every $k,l=1,2$.
\end{enumerate} 
\vspace{1em}  
Then there exist $z=(z_1,z_2)\in{\mathcal{H} }^2$ and a constant $C>0$, depending only on $x$, such that 
\vspace{1em}  
\begin{enumerate}
\item $x+z$ is an orthonormal system.
\vspace{1em}  
\item $z$ is an orthogonal system and every element of $z$ is orthogonal to every element of $x$.
\vspace{1em}  
\item $\norm{y_k-z_k}\leq C\sqrt{\epsilon} $ for every $k=1,2$.
\end{enumerate}  
\end{lemma} 
\subsection{Proof of the main result}
\begin{theorem}\label{exis}
For every $\delta\in[0,1]$ and for every $s\in(0,\infty)$, the supremum $S_{\delta}(s)$ is attained.
\end{theorem} 
\begin{proof}
The cases $\delta=0$ and $\delta=1$ have been proved in \cite{MR4480144} and \cite{nicola2026existenceextremizerssumeigenvalues}, respectively. Therefore, we focus only on the case $\delta\in(0,1)$ below.
\vspace{0.3em}  
\par By the definition of ${\mathcal{S} }_{\delta}(s)$, there exists a sequence $\{F^{(n)}: F^{(n)}=(f_1^{(n)},f_2^{(n)})\}_{n=1}^{\infty}\subset {\mathcal{F} }^2$ such that \[\lim\limits_{n\to\infty}{\mathcal{J} }_{\delta}(F^{(n)})={\mathcal{S} }_{\delta}(s),\]where each $F^{(n)}$ is an orthonormal system in ${\mathcal{F} }^2$. Obviously, for every $F\in{\mathcal{F} }^2$ we have \[{\mathcal{J} }_{\delta}(U_wF)={\mathcal{J} }_{\delta}(F),\quad w\in\mathbb{C}.\]Thus, by \autoref{fock lem}.(iv), we can assume that\[u_{F^{(n)}}(0)=\max\limits_{z\in\mathbb{C} }u_{F^{(n)}}(z),\quad n=1,2,\cdots.\]Since $\{F^{(n)}\}$ is a bounded sequence in ${\mathcal{F} }^2$, there exists $F=(f_1,f_2)\in{\mathcal{F} }^2$ such that $F^{(n)}\rightharpoonup F$. By the principle of uniform boundedness, we obtain $\norm{f_k}\leq1$ for every $k=1,2$. Furthermore, we have $F\neq0$. Otherwise, by the definition of $F$ and \autoref{fock lem}.(iii), we would have \[\lim\limits_{n\to\infty}u_{F^{(n)}}(0)=0\]and hence \[\lim\limits_{n\to\infty}\norm{u_{F^{(n)}}}_{\infty}=0,\]which implies that \[{\mathcal{S} }_{\delta}(s)=\limsup\limits_{n\to\infty}{\mathcal{J} }_{\delta}(F^{(n)})\leq\limsup\limits_{n\to\infty}s\norm{u_{F^{(n)}}}_{\infty}=0,\]a contradiction. According to \autoref{unitary}, we can suppose that $F$ is an orthogonal system and is ordered so that $\norm{f_2}\leq\norm{f_1}\leq1$.
\vspace{0.3em}  
\par Let $F^{(n)}-F=: G^{(n)}=(g_1,g_2)\in{\mathcal{F} }^2$, so that $G^{(n)}\rightharpoonup 0$. Notice that $F$ and $G^{(n)}$ satisfy assumptions in \autoref{hilbert}. By \autoref{lip} and these assumptions, we may assume that the elements of $G^{(n)}$ are pairwise orthogonal and orthogonal to all the elements of $F$.
\vspace{0.3em}  
\par Fix $\epsilon\in(0,1)$. For $n$ sufficiently large we have that 
\begin{equation}\label{supdef}
  {\mathcal{S} }_{\delta}(s)-\epsilon<{\mathcal{J} }_{\delta}(G^{(n)}+F)=\int_{{\Omega}_n}u_{G^{(n)}+F}(z) \ \dd A(z),
\end{equation}
where we let ${\Omega}_n={\Omega}_{G^{(n)}+F}(s)$. A direct computation yields that 
\begin{equation}\label{innerdecom}
u_{G^{(n)}+F}(z)=u_{G^{(n)}}(z)+u_F(z)+2\Re{\inner{F(z)}{G^{(n)}(z)}_{{\mathbb{C} }^2}}e^{-\pi|z|^2},\quad z\in\mathbb{C}.
\end{equation}
Here $\inner{\cdot}{\cdot}_{{\mathbb{C} }^2}$ denotes the $\delta$-weighted inner product on ${\mathbb{C} }^2$, i.e., \[\inner{\omega }{\nu}_{{\mathbb{C} }^2}={\omega}_1\bar{{\nu}}_1+\delta{\omega}_2\bar{{\nu}}_2,\quad \omega=({\omega}_1,{\omega}_2),\ \nu=({\nu}_1,{\nu}_2)\in{\mathbb{C} }^2.\]Since $u_F$ is integrable there exists a compact set $K\subset \mathbb{C} $ such that
\begin{equation}\label{Kdef}
\int_{\mathbb{C} \setminus K}u_F(z)\  \dd A(z)<\epsilon.  
\end{equation}
In addition, it follows from the Arzelá--Ascoli theorem that  weak convergence $G^{(n)}\rightharpoonup 0$ implies uniform convergence of $u_{G^{(n)}}$ on compact sets. Thus, for $n$ large enough it holds 
\begin{equation}\label{weak}
\int_{K}u_{G^{(n)}}(z)\  \dd A(z)<\epsilon.  
\end{equation}
Combining \eqref{Kdef} and \eqref{weak}, we obtain that
\[\begin{aligned}
  &\Bigl|\int_{{\Omega}_n}\Re{\inner{F(z)}{G^{(n)}(z)}_{{\mathbb{C} }^2}}e^{-\pi|z|^2} \ \dd A(z)\Bigr|\leq\int_{ {{\Omega}_n}\setminus K}\norm{F(z)}_{{\mathbb{C} }^2}\lVert G^{(n)}(z)\rVert_{\mathbb{C}^{2}}e^{-\pi|z|^2} \ \dd A(z)\\ &+\int_{ {{\Omega}_n}\cap K}\norm{F(z)}_{{\mathbb{C} }^2}\lVert G^{(n)}(z)\rVert_{\mathbb{C}^{2}}e^{-\pi|z|^2} \ \dd A(z)\leq\lVert G^{(n)}\rVert\Bigl(\int_{\mathbb{C} \setminus K}u_F(z)\  \dd A(z)\Bigr)^{1/2}+\\ &\norm{F}\Bigl(\int_{K}u_{G^{(n)}}(z)\  \dd A(z)\Bigr)^{1/2}\lesssim {\epsilon}^{\frac{1}{2}}. 
\end{aligned}\] 
Therefore, combining the estimate above with \eqref{supdef} and \eqref{innerdecom} we conclude 
\begin{equation}\label{1es}
{\mathcal{S} }_{\delta}(s)\leq {\mathcal{J} }_{\delta}(G^{(n)}+F)+\epsilon\leq \int_{{\Omega}_n}u_{G^{(n)}}(z) \ \dd A(z)+\int_{{\Omega}_n}u_{F}(z) \ \dd A(z)+C_{\delta}{\epsilon}^{\frac{1}{2}},
\end{equation}
where $C_{\delta}>0$ is a constant depending only on $\delta$. Using \autoref{decom} we observe 
\begin{equation}\label{dec}
{\mathcal{S} }_{\delta}(s)\leq \sum\limits_{k=0}^2a_k\int_{{\Omega}_n}u_{{\widehat{M}_k}^{(n)}}(z) \ \dd A(z)+C_{\delta}{\epsilon}^{\frac{1}{2}}.
\end{equation}
We claim that $f_2\neq0$. Otherwise, we obtain \[\begin{aligned}
  \sum\limits_{k=0}^2a_k\int_{{\Omega}_n}u_{{\widehat{M}_k}^{(n)}}(z) \ \dd A(z)&=(1-\norm{f_1}^2)\int_{{\Omega}_n}u_{{\widehat{G}}^{(n)}}(z) \ \dd A(z)+\norm{f_1}^2\int_{{\Omega}_n}u_{{\widehat{M}_1}^{(n)}}(z) \ \dd A(z)\\ &\leq (1-\norm{f_1}^2){\mathcal{S} }_{\delta}(s)+ \norm{f_1}^2\int_{{\Omega}_n}u_{{\widehat{M}_1}^{(n)}}(z) \ \dd A(z).
\end{aligned}\]
Since $F\neq0$, we have $\norm{f_1}>0$. From \eqref{dec} it follows 
\begin{equation}\label{2es}
{\mathcal{S} }_{\delta}(s)\leq \int_{{\Omega}_n}u_{{\widehat{M}_1}^{(n)}}(z) \ \dd A(z)+C_{\delta,F}{\epsilon}^{\frac{1}{2}}.
\end{equation}
Notice that $F$ is fixed and $\epsilon$ is arbitrary, \eqref{2es} implies that
\begin{equation}\label{conclu}
\lim\limits_{n\to\infty}\int_{{\Omega}_n}u_{{\widehat{M}_1}^{(n)}}(z) \ \dd A(z)={\mathcal{S} }_{\delta}(s).
\end{equation}
However, by the same arguments, for every $\varepsilon >0$ there exists a compact set $\widetilde{K}\subset \mathbb{C} $ such that 
\[\begin{aligned}
\int_{{\Omega}_n}u_{{\widehat{M}_1}^{(n)}}(z) \ \dd A(z)&\leq \int_{{\Omega}_n\cap \widetilde{K}}u_{\widehat{f_1}}(z) \ \dd A(z)+\delta\int_{{\Omega}_n\setminus \widetilde{K}}u_{\widehat{{g}_{2}^{(n)}}}(z) \ \dd A(z)+C_{\delta,F}{\varepsilon}^{\frac{1}{2}}\\ &\leq {\mathcal{S} }_0(|{\Omega}_n\cap \widetilde{K}|)+\delta{\mathcal{S} }_0(|{\Omega}_n\setminus \widetilde{K}|)+C_{\delta,F}{\varepsilon}^{\frac{1}{2}}
\end{aligned}\]
whenever $n$ is large enough. Thus, with \eqref{conclu} and \autoref{max} we obtain 
\[\begin{aligned}
{\mathcal{S} }_{\delta}(s)&\leq {\mathcal{S} }_0(|{\Omega}_n\cap \widetilde{K}|)+\delta{\mathcal{S} }_0(|{\Omega}_n\setminus \widetilde{K}|)\\ &\leq \max\limits_{t\in[0,s]}({\mathcal{S} }_0(t)+\delta{\mathcal{S} }_0(s-t))\\ &<{\mathcal{S} }_{\delta}(s),
\end{aligned}\]
a contradiction.
\vspace{0.3em}
\par Finally, arguing as in the derivation of \eqref{conclu}, we have 
\begin{equation}\label{conclu2}
\lim\limits_{n\to\infty}\int_{{\Omega}_n}u_{\widehat{F}}(z) \ \dd A(z)={\mathcal{S} }_{\delta}(s).
\end{equation}
Therefore, by the bathtub principle again, we conclude that $\widehat{F}$ is optimal, which completes the proof.
\end{proof}   
\section{The geometry of super--level sets}\label{geom}
In this section, we use $\Omega\subset \mathbb{C}$ to denote maximizer sets which achieves the supremum in \autoref{reduction}. The goal of this section is to give a geometric description of $\Omega$ for $\delta>0$ sufficiently small. In particular, we prove that, up to translation, such a set is a perturbation of a disc. The precise definition is as follows.
\begin{definition}
A set $\widetilde{\Omega}\subset \mathbb{C}$ is called a perturbation of a disc if there exist $R>0$ and $h\in{\mathcal{W} }_{\mathbb{R} }$ such that \[\widetilde{\Omega}={\Omega}_h :=\{re^{i\theta}:0\leq r<R+h(\theta), \theta\in[-\pi,\pi]\}.\] Here, ${\mathcal{W} }_{\mathbb{R} }$ denotes the Wiener algebra on the unit circle, i.e., \[{\mathcal{W} }_{\mathbb{R} }=\{f:\mathbb{R} \to\mathbb{R} :\sum\limits_{n\in\mathbb{Z} }|\widehat{f}(n)|<\infty\}.\]
\end{definition}
In addition, we will use the Fraenkel asymmetry of $\Omega$ to measure its distance to a disc:\[\mathcal{A}(\Omega) :=\inf \Bigl\{\frac{|\Omega\triangle B(z,r)|}{|\Omega|}: |B(x,r)|=|\Omega|, z\in\mathbb{C} \Bigr\}.\] \par In the remainder of this paper, we fix $R={\bigl(\frac{s}{\pi}\bigr)}^{\frac{1}{2}}$. Our main result reads as follows:
\begin{theorem}\label{geo}
There exists ${\delta}_s>0$ such that, for every $0\leq\delta\leq{\delta}_s$, one can find $z=z_{\delta,\Omega}\in\mathbb{C}$ and $h=h_{\delta,\Omega}\in{\mathcal{W} }_{\mathbb{R} }$ satisfying $\Omega-z={\Omega}_h$. Moreover, we have the following estimate \[\lim\limits_{\delta\to0}\norm{h}_{{\mathcal{W} }_{\mathbb{R} }}=0.\] 
\end{theorem}  
Before proving \autoref{geo}, let us begin with an easy observation.
\begin{lemma}\label{obs}
For every $0\leq\delta\leq1$, we have
\begin{equation}\label{eqobs}
0\leq {\lambda }_1(s)-{\lambda }_1(\Omega)\lesssim_s\delta. 
\end{equation}
Here, for every $j=1,2,\cdots$, ${\lambda }_j(s)$ is defined as \[{\lambda }_j(s):=1-e^{-s}\sum\limits_{n=0}^{j-1}\frac{s^n}{n!}.\]
\end{lemma}  
\begin{proof}
We claim that ${\lambda}_1(s)\leq(1+\delta){\lambda }_1(\Omega)$. Otherwise, we have \[{\lambda }_1(\Omega)+\delta{\lambda }_2(\Omega)\leq(1+\delta){\lambda }_1(\Omega)<{\lambda }_1(s)<{\mathcal{S} }_{\delta}(s),\]which is contrary to the choice of $\Omega$. Thus, by \cite[Theorem 3.1]{MR4480144}, we obtain \[0\leq{\lambda }_1(s)-{\lambda }_1(\Omega)\leq\delta{\lambda }_1(\Omega)\leq\delta{\lambda }_1(s),\]which yields the desired result.
\end{proof}  
According to the stability of the Faber--Krahn inequality for STFT \cite[Theorem 1.5]{MR4728242}, we have\[\lim\limits_{\delta\to0}\mathcal{A}(\Omega)=0,\]which implies that for every $\varepsilon >0$ there exists $z_{\varepsilon }\in\mathbb{C} $ such that $|\Omega\triangle B(z_{\varepsilon },R)|<\varepsilon$. Notice that eigenvalues are invariant under translations, i.e., \[{\lambda}_j(\Omega)={\lambda}_j(\Omega-w),\quad w\in\mathbb{C} \ \text{and} \ j=1,2,\cdots.\]Therefore, to prove \autoref{geo}, it suffices to prove the following theorem:
\begin{theorem}\label{geo1}
Suppose that $\Omega$ is a maximizer set at which the supremum in \autoref{reduction} is attained and satisfies 
\begin{equation}\label{sym}
\lim\limits_{\delta\to0}|\Omega\triangle B(0,R)|=0.
\end{equation}
Then there exists ${\delta}_s>0$ such that $\delta<{\delta}_s$ implies the existence of $h=h_{\delta,\Omega}\in{\mathcal{W} }_{\mathbb{R} }$ with the property that \[\Omega={\Omega}_h\]and
\begin{equation}\label{h}
\lim\limits_{\delta\to0}\norm{h}_{{\mathcal{W} }_{\mathbb{R} }}=0.
\end{equation}
\end{theorem} 
\subsection{Preliminary estimates} 
For simplicity, we use $T_0$ to denote the Toeplitz oparator $T_{B(0,R)}$. In addition, we warn the reader that most of the functions and variables in this section depend on $\delta$. However, we suppress this dependence from the notation. 
\vspace{0.3em}
\par It is well known \cite{MR966733, MR1091688} that the eigenfunctions of $T_0$ are $\{e_n\}_{n=0}^{\infty}$ and the corresponding
eigenvalues are $\{{\lambda}_n(s)\}_{n=1}^{\infty}$, respectively. Thus, every eigenvalue of $T_0$ is simple. By \eqref{s1}, we obtain 
\begin{equation}\label{egv}
  |{\lambda}_n(\Omega)-{\lambda}_n(s)|\leq\norm{T_{\Omega}-T_0}\leq|\Omega\triangle B(0,R)|,\quad n=1,2,\cdots,
  \end{equation}
  which reveals that every eigenvalue of $T_{\Omega}$ is also simple. Let $f,g\in{\mathcal{F}}$ be the first two unit eigenfunctions of $T_{\Omega}$ corresponding to ${\lambda}_1(\Omega)$ and ${\lambda}_2(\Omega)$, respectively. According to \autoref{reduction}, we have $\Omega={\Omega}_F(s)$, where $F=(f,g)\in{\mathcal{F} }^2$. Moreover, we have the following result: 
\begin{lemma}\label{egf}
Under the same assumptions as in \autoref{geo1}, it holds 
\begin{equation}\label{eqegf}
\inf\limits_{c\in\mathbb{C}:|c|=1}\norm{f-c}\lesssim_s\delta.
\end{equation}  
\begin{proof}
Notice that \[\inf\limits_{c\in\mathbb{C}:|c|=1}\norm{f-c}^2=2(1-|\inner{f}{1}|)\leq2(1-|\inner{f}{1}|^2)=2\norm{\mathcal{P}f}^2,\]where $\mathcal{P} $ denotes the projection from $\mathcal{F} $ onto $\overline{\text{span}}\{e_n:n=1,2,\cdots\}$. Since ${\lambda}_1(s)$ is simple, we have \[m(s) \ :=\inf\limits_{n\geq1}({\lambda}_1(s)-{\lambda}_n(s))>0\]and hence by \eqref{egv} we obtain \[\widetilde{m}(s) \ :=\inf\limits_{n\geq1}({\lambda}_1(\Omega)-{\lambda}_n(s))>0\]provided $\delta$ is sufficiently small. In addition, we obtain that \[{\widetilde{m}(s)}^2\norm{\mathcal{P}f}^2\leq\norm{\mathcal{P} (T_{\Omega}-T_0)f}^2\lesssim_s{\delta}^2,\]which, combined with the previous estimate , yields the desired conclusion.
\end{proof}   
\begin{remark}
It is easy to see $u_{\alpha G}=u_G$ for every $|\alpha |=1$ and $G\in{\mathcal{F} }^2$. Hence, without loss of generality, we may assume  
\begin{equation}\label{eqegf'}
\norm{f-e_0}\lesssim_s\delta.
\end{equation}  
\end{remark}   
\end{lemma}
\begin{corollary}\label{uFes}
Under the same assumptions as in \autoref{geo1}, we have 
\begin{equation}\label{equF}
\lim\limits_{\delta\to0}\norm{u_F-\mathfrak{u} }_{\infty}=0.
\end{equation}  
Here, $\mathfrak{u}(z)=e^{-\pi|z|^2} $ for each $z\in\mathbb{C} $.
\end{corollary}  
\begin{proof}
By \autoref{fock lem}.(iii) we obtain \[\norm{u_F-\mathfrak{u}}_{\infty}\leq\norm{u_f-\mathfrak{u}}_{\infty}+\delta.\]Notice that for every $z\in\mathbb{C} $ it holds 
\[\begin{aligned}
  |u_f(z)-\mathfrak{u}(z)|&\leq(|f(z)|+1)|f(z)-1|e^{-\pi|z|^2}\\ &=(|\inner{f}{k_z}|+e^{-\pi|z|^2})|\inner{f-e_0}{k_z}|\\ &\lesssim \norm{f-e_0}.
\end{aligned}\]  
which, combined with \eqref{eqegf'}, concludes the proof.
\end{proof}  
Let $t>0$ be the unique number such that $\Omega=\{u_F>t\}$. It is to see that $B(0,R)=\{\mathfrak{u}>e^{-s}\}$. Our next goal is to bound the difference between $t$ and $e^{-s}$. Indeed, we have the following estimate:
\begin{lemma}\label{df}
Under the same assumptions as in \autoref{geo1}, it follows that 
\begin{equation}\label{eqdf}
|t-e^{-s}|\leq5\norm{u_F-\mathfrak{u} }_{\infty}.
\end{equation}
\end{lemma}  
\begin{proof}
If $|t-e^{-s}|>5\norm{u_F-\mathfrak{u} }_{\infty}$, then there are only two possible cases:
\vspace{0.5em}
\par \textbf{Case 1:} $t>e^{-s}+5\norm{u_F-\mathfrak{u} }_{\infty}$. In this case, we have \[\Omega=\{u_F>t\}\subset \{\mathfrak{u}>e^{-s}+4\norm{u_F-\mathfrak{u} }_{\infty}\},\]which implies \[|\Omega|\leq|\{\mathfrak{u}>e^{-s}+4\norm{u_F-\mathfrak{u} }_{\infty}\}|<s,\]a contradiction.
\vspace{0.5em}
\par \textbf{Case 2:} $e^{-s}>t+5\norm{u_F-\mathfrak{u} }_{\infty}$. In this case, we have \[\{\mathfrak{u}>e^{-s}-4\norm{u_F-\mathfrak{u} }_{\infty}\}\subset \{u_F>t\}=\Omega,\]which implies that \[s<|\{\mathfrak{u}>e^{-s}-4\norm{u_F-\mathfrak{u} }_{\infty}\}|\leq|\Omega|,\]a contradiction.
\end{proof}  
\begin{corollary}\label{df'}
 Under the same assumptions as in \autoref{geo1}, we have
\begin{equation}\label{eqdf'}
B(0,R/2)\subset \Omega\subset B(0,2R)
\end{equation}
provided $\delta$ is sufficiently small. 
\end{corollary}  
\begin{proof}
For every $z\in B(0,R/2)$, combining \autoref{uFes} and \autoref{df}, we obtain
\[\begin{aligned}
  u_F(z)&\geq\mathfrak{u}(z)-\norm{u_F-\mathfrak{u} }_{\infty}\\ &\geq e^{-s/4}-\norm{u_F-\mathfrak{u} }_{\infty}\\ & \geq t+ (e^{-s/4}-e^{-s}-6\norm{u_F-\mathfrak{u} }_{\infty})\\ &>t
\end{aligned}\]     
whenever $\delta$ is sufficiently small. In addition, if $|z|\geq 2R$, then we have 
\[\begin{aligned}
  u_F(z)&\leq\mathfrak{u}(z)+\norm{u_F-\mathfrak{u} }_{\infty}\\ &\leq e^{-4s}+\norm{u_F-\mathfrak{u} }_{\infty}\\ & \leq t+ (e^{-4s}-e^{-s}+6\norm{u_F-\mathfrak{u} }_{\infty})\\ &<t
\end{aligned}\]
if $\delta$ is sufficiently small.
\end{proof}  
\subsection{Proof of \autoref{geo1}}
In this section, we split the proof into two steps.
\vspace{0.3em}
\par \textbf{Step 1.} Given $\theta \in[-\pi,\pi)$, we define 
\[\begin{aligned}{\mathcal{G}_{\theta} }(r)\ &:=u_F(re^{i\theta}),\\{\mathcal{G}_{\theta1} }(r)\ &:=u_f(re^{i\theta}),\\{\mathcal{G}_{\theta2} }(r)\ &:=u_g(re^{i\theta})
\end{aligned}\]
for each $r\in[0,\infty)$. By the Cauchy--Schwarz inequality, for every $\tilde{f}\in\mathcal{F} $ with $\lVert \tilde{f} \rVert=1$, we have 
\begin{equation}\label{fkes}
  |{\tilde{f}}'(z)|\lesssim_s\lVert \tilde{f} \rVert,\quad |z|\leq 10R.
\end{equation}  
Thus, according to \autoref{fock lem}.(iii) and \eqref{fkes}, it follows that 
\begin{equation}\label{g2es}
  |{\mathcal{G}_{\theta2} }'(r)|\leq2\bigl(|g'(re^{i\theta})g(re^{i\theta})|+\pi r\bigr)e^{-\pi r^2}\lesssim_s\norm{g}=1
\end{equation}  
provided $r\in[0,2R]$. Similarly, thanks to \eqref{eqegf'}, we obtain 
\begin{equation}\label{g1es}
\begin{aligned}
  {\mathcal{G}_{\theta1} }'(r)&\leq2\bigl(|f'(re^{i\theta})f(re^{i\theta})|-\pi r\bigr)e^{-\pi r^2}\\ &=\bigl(|(f-e_0)'(re^{i\theta})f(re^{i\theta})|-\pi r\bigr)e^{-\pi r^2}\\ &\leq-\pi re^{-\pi r^2}<0
\end{aligned}  
\end{equation} 
whenever $r\in[R/2,2R]$ and $\delta$ is small enough. Combining \eqref{eqdf'}, \eqref{g2es} and \eqref{g1es} it holds 
\begin{equation}\label{ges}
  |{\mathcal{G}_{\theta} }'(r)|\lesssim_sRe^{-\pi R^2}<0
\end{equation}  
for every $z=re^{i\theta}\in\Omega$ and $\delta$ close to 0. Standard arguments yield that ${\mathcal{G}_{\theta} }(\frac{R}{2})<t$ and ${\mathcal{G}_{\theta} }(2R)>t$. Thus, by \eqref{ges} there is a unique $r(\theta)\in[R/2,2R]$ such that ${\mathcal{G}_{\theta} }(r(\theta))=t$. Wrriting $h(\theta) :=r(\theta)-R$, we have therefore proved the existence and uniqueness of $h$. Furthermore, by the same arguments, we obtain that there exists a function $\eta: [0,\infty)\rightarrow  [0,\infty)$ satisfying \[\lim\limits_{\delta\to0}\eta(\delta)=0\] and \[r(\theta)\in[(1-\eta(\delta))R,(1+\eta(\delta))R].\]In particular, for every $r\in[R/2,2R]$, by \eqref{g2es} and \eqref{g1es}, there exist $\delta_{s}>0$ and $K_s>0$ such that $\delta\in[0,\delta_s]$ implies $|{\mathcal{G}_{\theta} }'(r)|\geq K_s$. Let $\eta(\delta) :=10{K_s}^{-1}\norm{u_F-\mathfrak{u} }_{\infty}$. According to \eqref{equF} and \eqref{eqdf}, if $|r-R|\geq\eta(\delta)R$, then we have 
\[\begin{aligned}
|u_F(re^{i\theta})-t|&\geq|u_F(re^{i\theta})-u_F(Re^{i\theta})|-|u_F(Re^{i\theta})-\mathfrak{u} (Re^{i\theta})|-|\mathfrak{u} (Re^{i\theta})-t|\\ &\geq|\mathcal{G}_{\theta}(r)-\mathcal{G}_{\theta}(R)|-\norm{u_F-\mathfrak{u}}_{\infty}-|e^{-s}-t|\\ &\geq 4\norm{u_F-\mathfrak{u}}_{\infty}>0,
\end{aligned}\]  
a contradiction.
\vspace{0.3em}
\par \textbf{Step 2.} For every $\theta\in[-\pi,\pi)$ and $r\in[0,\infty)$, we define \[\mathfrak{G} (r,\theta) :=u_F(re^{i\theta}).\]Obviously, $\mathfrak{G} $ is a smooth function on $[0,\infty) \times \mathbb{T} $ and $\mathfrak{G} (r(\theta),\theta)=t$. By the computations above, we have $\partial_r \mathfrak{G} (r,\theta)<0 $ provided $\delta$ is small enough. Thus, according to the smooth implicit function theorem \cite[Theorem C.40]{lee2013introduction}, there exist neighborhoods $V_{\theta}\subset \mathbb{T} $ of $\theta$ and $I_{\theta}\subset [0,\infty)$ of $r(\theta)$ and a smooth function $\tilde{r}: I_{\theta}\rightarrow V_{\theta}$ such that $\mathfrak{G} (\tilde{r}(\vartheta ),\vartheta )=t$ for every $\theta\in V_{\theta}$. Let $V_{\theta}$ sufficiently small such that\[|\tilde{r}(\vartheta )-r(\theta)|<\eta(\delta)R,\quad \vartheta \in V_{\theta},\] and hence $\tilde{r}(\vartheta )\in[R/2,2R]$ for every $\vartheta \in V_{\theta}$. By the uniqueness of $r$, we conclude that $\tilde{r}=r$ in $V_{\theta}$, which implies that $h\in C^{\infty}(\mathbb{T} )$ and 
\begin{equation}\label{hes}
\norm{h}_{\infty}\leq\eta(\delta)R.
\end{equation}
Notice that \[h'(\theta)=r'(\theta)=-\frac{\partial_{\theta} \mathfrak{G} (r,\theta)}{\partial_r \mathfrak{G} (r,\theta)}=-\frac{\partial_{\theta} \mathfrak{G} (r,\theta)}{\mathcal{G}_{\theta}'(r)}.\]It follows by the same arguments as in \textbf{Step 1} that $\partial_{\theta} \mathfrak{G} (r,\theta)$ tends to 0 uniformly on $[R/2,2R]\times\mathbb{T} $ as $\delta\to0$. Thus, we have 
\begin{equation}\label{h'es}
\lim\limits_{\delta\to0}\norm{h'}_{\infty}=0.
\end{equation}
Combining the Cauchy--Schwarz inequality and the Parseval's identity, it holds \[\begin{aligned}
\sum\limits_{n\in\mathbb{Z} }|\widehat{h}(n)|&=|\widehat{h}(0)|+\sum\limits_{|n|=1}^{\infty}\frac{1}{2\pi n}|2\pi n\widehat{h}(n)|\\ &=|\widehat{h}(0)|+\sum\limits_{|n|=1}^{\infty}\frac{1}{2\pi |n|}|\widehat{h'}(n)|\\  &\leq\norm{h}_{\infty}+\Bigl(\sum\limits_{|n|=1}^{\infty}\frac{1}{4{\pi}^2 n^2}\Bigr)^{\frac{1}{2}}\Bigl(\sum\limits_{|n|=1}^{\infty}|\widehat{h'}(n)|^2\Bigr)^{\frac{1}{2}}\\ &\lesssim\norm{h}_{\infty}+\norm{h'}_2\\ &\leq\norm{h}_{\infty}+\norm{h'}_{\infty},
\end{aligned}\]
which, combined with \eqref{hes} and \eqref{h'es}, concludes the proof.
\section{Analyticity of vector-valued mappings}\label{analy}
In this section, we will prove \autoref{main}. Our strategy is to use the vector-valued version of the analytic implicit function theorem \cite[Theorem 15.3]{MR787404} to conclude that the perturbation function $h$ obtained in \autoref{geo} vanishes whenever $\delta>0$ is small enough. In what follows, we continue to use the same notations as in the previous sections.
\subsection{Analyticity of \texorpdfstring{$T_{{\Omega}_h}$}{TOmegah} and of the corresponding eigenvalues and eigenfunctions} 
 For any $x,y$ in a Hilbert space $\mathcal{H}$, we define the tensor product of $x$ and $y$ as\[(x\otimes y )(z) := \inner{z}{y}x,\quad z\in\mathcal{H}.\]Obviously, we have $x\otimes y\in\mathfrak{S}_1(\mathcal{H} )$ and $\norm{x\otimes y}_{\mathfrak{S} _1}\leq\norm{x}\norm{y}$. For every $w\in\mathbb{C} $, we use $\mathcal{Q}_w$ to denote $k_w\otimes k_w$. It is easy to check that the mapping $w\mapsto \mathcal{Q}_w$ is continuous from $\mathbb{C}$ to $\mathfrak{S}_1$. Thus, the Bochner integral $\int_{\Omega}\mathcal{Q}_w\ \dd A(w)$ is well defined \cite{MR3617205} and in fact coincides with $T_{\Omega}$. Notice that 
\begin{equation}\label{hdf}
\lim\limits_{\norm{h}_{\infty}\to0}|\Omega_{h}\triangle B(0,R)|=0,
\end{equation}  
which implies  
\begin{equation}\label{hdf'}
\lim\limits_{\norm{h}_{\infty}\to0}\norm{T_{\Omega_h}-T_{B(0,R)}}=0.
\end{equation} 
 We now show the following result:
\begin{proposition}\label{anaop}
 There exists $\varepsilon>0$ such that the mapping \[\mathcal{T}: {\mathcal{W} }_{\mathbb{R} }\rightarrow \mathcal{L} (\mathcal{F} ), \ h\mapsto T_{\Omega_h}\]is analytic whenever $\norm{h}_{{\mathcal{W} }_{\mathbb{R} }}<\varepsilon $. Here, $\mathcal{L}(\mathcal{H} )$ denotes the space of bounded linear operators on $\mathcal{H} $. Furthermore, we have 
 \begin{equation}\label{eqop}
  \mathcal{D} \mathcal{T} (0)[h]=R\int_{-\pi}^{\pi}h(\theta)\mathcal{Q}_{Re^{i\theta}}\ \dd \theta.
 \end{equation} 
\end{proposition}  
\begin{proof}
Let $\varepsilon=\frac{1}{2}\min\{R,1\}$. Formally, by polar-coordinate decomposition formula, we have 
\[\begin{aligned}
\mathcal{T}(h)&=\int_{\Omega}\mathcal{Q}_w\ \dd A(w)\\ &=\int_{-\pi}^{\pi}\int_{0}^{R+h(\theta)}r\mathcal{Q}_{re^{i\theta}}\ \dd r\ \dd \theta\\ &=\int_{-\pi}^{\pi}\int_{0}^{R+h(\theta)}r\Bigl[e^{-\pi r^2/2}\sum\limits_{m=0}^{\infty}\Bigl(\frac{{\pi}^m}{m!}\Bigr)^{1/2}r^me^{-im\theta}e_m\Bigr] \\ &\otimes\Bigl[e^{-\pi r^2/2}\sum\limits_{n=0}^{\infty}\Bigl(\frac{{\pi}^n}{n!}\Bigr)^{1/2}r^ne^{-in\theta}e_n\Bigr]\ \dd r\ \dd \theta\\ &=\sum\limits_{m,n=0}^{\infty}\mathcal{B}_{m,n}(h)(e_m\otimes e_n),
\end{aligned}\]
where \[\mathcal{B}_{m,n}(h)=\Bigl(\frac{\pi^{m+n}}{m!n!}\Bigr)^{1/2}\int_{-\pi}^{\pi}\int_{0}^{R+h(\theta)}e^{-\pi r^2}r^{m+n+1}e^{-i(m+n)\theta}\ \dd r\dd\theta.\]Therefore, to prove the analyticity of $\mathcal{T} $, it suffices to to show that the above infinite series is absolutely summable whenever $\norm{h}_{{\mathcal{W} }_{\mathbb{R} }}<\varepsilon $. Indeed, it holds \[\sum\limits_{m,n=0}^{\infty}|\mathcal{B}_{m,n}(h)|\lesssim\sum\limits_{m,n=0}^{\infty}\Bigl(\frac{\pi^{m+n}}{m!n!}\Bigr)^{1/2}\int_{0}^{R+1}e^{-\pi r^2}r^{m+n+1}\ \dd r\lesssim_s1,\]which yields both the desired conclusion and \eqref{eqop} by routine computations.
\end{proof}  
We are now ready to prove the analyticity of eigenfunctions and eigenvalues.
\begin{proposition}\label{anaei}
For each $j=1,2,\cdots$, there exist $\delta_j>0$ and analytic mappings \[f_j:B_{{\mathcal{W} }_{\mathbb{R} }}(0,\delta_j)\rightarrow \mathcal{F}\]and \[\lambda_j:B_{{\mathcal{W} }_{\mathbb{R} }}(0,\delta_j)\rightarrow \mathbb{C} \] such that $f_j(h)$ and $\lambda_j(h)$ are the $j$-th unit eigenfunction and the $j$-th eigenvalue of $T_{{\Omega}_h}$, respectively, satisfying $f_j(0)=e_{j-1}$ and $\lambda_j(0)=\lambda_j(B(0,R))$.
\end{proposition} 
\begin{proof}
Let $\varepsilon $ be the constant in \autoref{anaop} and \[\tilde{\delta}_j=\frac{1}{K_s}\min\{\varepsilon ,\lambda_{j-1}(s)-\lambda_j(s), \lambda_{j}(s)-\lambda_{j+1}(s)\}.\] Here $K_s$ is a sufficiently large constant to be chosen later. For every $\lambda\in\mathbb{C} $, $f\in\mathcal{F}\ominus\mathbb{C} e_{j-1} $, and $h\in{\mathcal{W} }_{\mathbb{R} }$,  we define\[\Phi (\lambda,f,h) :=\mathcal{T} (h)(e_{j-1}+f)-\lambda(e_{j-1}+f).\]Thanks to \autoref{anaop}, we obtain that $\Phi$ is an analytic mapping from $\mathbb{C}\times(\mathcal{F}\ominus\mathbb{C}e_{j-1})\times{\mathcal{W} }_{\mathbb{R} }$ to $\mathcal{F} $ whenever $\norm{h}_{{\mathcal{W} }_{\mathbb{R} }}<\tilde{\delta}_j$ and $\Phi(\lambda_{j}(s),0,0)=0$. Furthermore, it holds 
\begin{equation}\label{eqdrv}
\partial_{(\lambda,f)}\Phi(\lambda_{j}(s),0,0)[\eta,\nu]= (T_0-\lambda_{j}(s))\nu-\eta e_j,
\end{equation}
where $T_0=\mathcal{T}(0)=T_{B(0,R)}$. Notice that for every $\nu\in(\mathcal{F}\ominus\mathbb{C}e_{j-1})$, we have $\inner{\nu}{e_{j-1}}=0$, which implies that \[(T_0-\lambda_{j}(s))\nu=\sum\limits_{n\neq j}(\lambda_n(s)-\lambda_{j}(s))\inner{\nu}{e_n}e_n.\]Since the eigenvalue $\lambda_{j}(s)$ is simple, we obtain that \[\inf\limits_{n\neq j}|\lambda_n(s)-\lambda_{j}(s)|>0\] and hence \[(T_0-\lambda_{j}(s))^{-1}\nu :=\sum\limits_{n\neq j}(\lambda_n(s)-\lambda_{j}(s))^{-1}\inner{\nu}{e_n}e_n\in(\mathcal{F}\ominus\mathbb{C}e_{j-1}),\]which yields that $\partial_{(\lambda,f)}\Phi(\lambda_{j}(s),0,0)$ is invertible from $\mathbb{C}\times(\mathcal{F}\ominus\mathbb{C}e_{j-1})$ to $\mathcal{F} $ and in particular \[\Bigl(\partial_{(\lambda,f)}\Phi(\lambda_j(s),0,0)\Bigr)^{-1}(g)=\Bigl[-\inner{g}{e_{j-1}},(T_0-\lambda_j(s))^{-1}(g-\inner{g}{e_{j-1}}e_{j-1})\Bigr],\quad g\in\mathcal{F}.\]Therefore, according to the vector-valued version of the analytic implicit function theorem \cite[Theorem 15.3]{MR787404}, there exist $\hat{\delta}_j>0$ and unique analytic mappings \[{\tilde{f}}_j:B_{{\mathcal{W} }_{\mathbb{R} }}(0,\hat{\delta}_j)\rightarrow \mathcal{F}\]and \[\lambda_j:B_{{\mathcal{W} }_{\mathbb{R} }}(0,\hat{\delta}_j)\rightarrow \mathbb{C} \] such that $\Phi(\lambda_j(h),{\tilde{f}}_j(h),h)=0$ for every $h\in B_{{\mathcal{W} }_{\mathbb{R} }}(0,\hat{\delta}_j)$. Let \[\delta_j=\min\{\tilde{\delta}_j,\hat{\delta}_j\}\] and \[f_j(h)=\frac{\tilde{f}_j+e_{j-1}}{\norm{\tilde{f}_j+e_{j-1}}}=\frac{\tilde{f}_j+e_{j-1}}{\sqrt{1+\norm{\tilde{f}_j(h)}^2} }.\]It is easy to see that $f_j$ is analytic and satisfies \[\mathcal{T}(h)[f_j(h)]=\lambda_j(f)f_j(h),\quad \norm{h}_{{\mathcal{W} }_{\mathbb{R} }}<\delta_j.\]

Finally, it remains to show that $\lambda_j(h)$ is the $j$-th eigenvalue of $\mathcal{T}(h)$. Indeed, according to \eqref{egv}, \eqref{hdf'}, and the definition of $\delta_j$, it follows that all the eigenvalues of $\mathcal{T} (h)$ must be simple whenever $K_s$ is large enough and $\norm{h}_{{\mathcal{W} }_{\mathbb{R} }}<\delta_j$. In addition, denoting by $\{\lambda_n(\Omega_h)\}_{n=1}^{\infty}$ all the eigenvalues of $\mathcal{T} (h)$, it holds 
\[\begin{aligned}
\inf\limits_{n\neq j}|\lambda_n(\Omega_h)-\lambda_j(s)|&\geq\inf\limits_{n\neq j}|\lambda_n(\Omega_h)-\lambda_j(\Omega_h)|-|\lambda_j(\Omega_h)-\lambda_j(s)|\\ &\geq\inf\limits_{n\neq j}|\lambda_n(s)-\lambda_j(s)|-5\max\limits_{n\geq1}|\lambda_n(\Omega_h)-\lambda_n(s)|\\ &\geq\frac{1}{2}\min\bigl\{\lambda_{j-1}(s)-\lambda_j(s),\lambda_j(s)-\lambda_{j+1}(s)\bigr\}>0.
\end{aligned}\]
Notice that $\lambda_j(h)\in\{\lambda_n(\Omega_h)\}_{n=1}^{\infty}$. By the continuity of $\lambda_j$ with respect to $h$, we have \[|\lambda_j(h)-\lambda_j(s)|<\inf\limits_{n\neq j}|\lambda_n(\Omega_h)-\lambda_j(s)|\] whenever $\hat{\delta}_j$ is chosen sufficiently small and $\norm{h}_{{\mathcal{W} }_{\mathbb{R} }}<\hat{\delta}_j$. Thus, we have $\lambda_j(h)=\lambda_j(\Omega)$, which concludes the proof.
\end{proof} 
\begin{corollary}\label{inner}
For every $h\in B_{{\mathcal{W} }_{\mathbb{R} }}(0,\delta_j)$, we have 
\begin{equation}\label{eqin}
\inner{\mathcal{D}f_j(0)[h]}{e_{j-1}}=0.
\end{equation}
\end{corollary}  
\begin{proof}
 Recalling the construction of $f_j$, we obtain that \[\inner{f_j(\mu h)[\phi]}{e_{j-1}}>0\] for every $\mu \in\mathbb{R}$ and \[\inner{f_j(0)[\phi]}{e_{j-1}}=0.\]Thus, the desired conclusion follows from Fermat's lemma.
\end{proof}  
\subsection{Further geometric properties of the maximizer set and related analytic properties}
Review that $\Omega$ is the maximizer set achieving the supremum in \autoref{reduction}. According to the arguments in Section 4, we may assume without loss of generality that \eqref{sym} holds. Additionally, we can suppose that 
\begin{equation}\label{ave}
z_{\Omega} :=\int_{\Omega}z\ \dd A(z)=0.
\end{equation} 
Otherwise, we may use replace $\Omega$ by $\Omega-z_{\Omega}$. Indeed, notice that eigenvalues are invariant under translations, and \[\Bigl|\int_{\Omega}z\ \dd A(z)\Bigr|=\Bigl|\int_{\Omega}z\ \dd A(z)-\int_{B(0,R)}z\ \dd A(z)\Bigr|\leq\int_{\Omega\triangle B(0,R)}|z|\ \dd A(z).\] Thus, by \autoref{df'}, $\Omega-z_{\Omega}$ also satisfies \eqref{sym}. According to \autoref{geo1}, there exists $h_{\delta}\in{\mathcal{W} }_{\mathbb{R} }$ such that $\Omega=\Omega_{h_{\delta}}$. Since $h_{\delta}$ is a real-valued function, there exist $\alpha_{\delta},\beta_{\delta},\gamma_{\delta}  \in\mathbb{R} $ and $\varphi_{\delta} \in{\mathcal{W} }_{\mathbb{R} }$ satisfying
\begin{equation}\label{van}
\widehat{\varphi_{\delta} }(0)=\widehat{\varphi_{\delta}}(1)=\widehat{\varphi_{\delta} }(-1)=0.
\end{equation}
and $h(\theta)=\alpha_{\delta}+\beta_{\delta}\cos\theta+\gamma_{\delta} \sin\theta+\varphi_{\delta} (\theta)$. Combining \eqref{ave}, $|\Omega_{h_{\delta}}|=s$, and polar-coordinate decomposition formula, we have 
\begin{equation}\label{eqgeo2}
  \begin{aligned}
  &\frac{1}{2}\int_{-\pi}^{\pi}\bigl(R+\alpha_{\delta}+\beta_{\delta}\cos\theta+\gamma_{\delta} \sin\theta+\varphi_{\delta} (\theta)\bigr)^2\ \dd \theta=s,\\ &\frac{1}{3}\int_{-\pi}^{\pi}\bigl(R+\alpha_{\delta}+\beta_{\delta}\cos\theta+\gamma_{\delta} \sin\theta+\varphi_{\delta} (\theta)\bigr)^3\cos\theta\ \dd \theta=0,\\ &\frac{1}{3}\int_{-\pi}^{\pi}\bigl(R+\alpha_{\delta}+\beta_{\delta}\cos\theta+\gamma_{\delta} \sin\theta+\varphi_{\delta} (\theta)\bigr)^3\sin\theta\ \dd \theta=0.
\end{aligned}
\end{equation}
We are now ready to introduce the following result:
\begin{proposition}\label{geo2}
Set $\mathcal{W}_0 :=\{\varphi\in{\mathcal{W} }_{\mathbb{R} }: \varphi \ \text{satisfies} \ \eqref{van}\}$. Then there exist $\iota>0$ and analytic mappings $\alpha$, $\beta$, and $\gamma$, where
\[
\alpha,\beta,\gamma:B_{\mathcal{W}_0}(0,\iota)\to\mathbb{R},
\]
satisfying \eqref{eqgeo2}. Furthermore, it holds that
\begin{equation}\label{eqvan}
\alpha(0)=\beta(0)=\gamma(0)=0
\end{equation}
and
\begin{equation}\label{eqdrvan}
(\mathcal{D}\alpha)(0)=(\mathcal{D}\beta)(0)=(\mathcal{D}\gamma)(0)=0.
\end{equation}
\end{proposition}  
\begin{proof}
For every $(\alpha,\beta,\gamma,\varphi)\in{\mathbb{R} }^3\times{\mathcal{W} }_0$, set
\[
\Psi(\alpha,\beta,\gamma,\varphi)=\begin{pmatrix}
  \frac{1}{2}\int_{-\pi}^{\pi}\bigl(R+h(\alpha,\beta,\gamma,\varphi)(\theta)\bigr)^2\ \dd \theta \\[3pt]
  \frac{1}{3}\int_{-\pi}^{\pi}\bigl(R+h(\alpha,\beta,\gamma,\varphi)(\theta)\bigr)^3\cos\theta\ \dd \theta \\[3pt]
  \frac{1}{3}\int_{-\pi}^{\pi}\bigl(R+h(\alpha,\beta,\gamma,\varphi)(\theta)\bigr)^3\sin\theta\ \dd \theta
\end{pmatrix},
\]
where $h(\alpha,\beta,\gamma,\varphi)(\theta)=\alpha+\beta\cos\theta+\gamma \sin\theta+\varphi (\theta)$. By standard computations, we have 
\[
\Psi(0,0,0,0)=\begin{pmatrix}
  s \\
  0 \\
  0
\end{pmatrix} 
\]
and 
\[
J_{(\alpha,\beta,\gamma)}\Psi(0,0,0,0)= \begin{pmatrix}
  2\pi R & 0 & 0 \\
  0 & \pi R^2 & 0 \\
  0 & 0 & \pi R^2
\end{pmatrix}.
\] 
Thus, by the analytic implicit function theorem, we complete the proofs of analyticity and \eqref{eqvan}. 
\vspace{0.3em}
\par Next, we shall compute $(\mathcal{D}\alpha)(0)$, and the remaining arguments are analogous. For all $\mu \in\mathbb{R} $ and $\varphi\in\mathcal{W}_0$, we have
\[
\frac{1}{2}\int_{-\pi}^{\pi}\bigl(R+\alpha(\mu \varphi)+\beta(\mu \varphi)\cos\theta+\gamma(\mu \varphi) \sin\theta+\mu \varphi (\theta)\bigr)^2\ \dd \theta=s,
\]
which implies that 
\[
\int_{-\pi}^{\pi}R\bigl[(\mathcal{D}\alpha)(0)[\varphi]+(\mathcal{D}\beta)(0)[\varphi]\cos\theta+(\mathcal{D}\gamma)(0)[\varphi]\sin\theta+\varphi(\theta)\bigr]\ \dd\theta=0.
\]
Notice that $\int_{-\pi}^{\pi}\varphi(\theta)\ \dd\theta=0$ since $\varphi\in\mathcal{W}_0$. Thus, we obtain $(\mathcal{D}\alpha)(0)[\varphi]=0$, which yields the desired conclusion. 
\end{proof}  
\begin{remark}\label{hpro}
It is easy to see that $h(0)=0$ and $(\mathcal{D}h)(0)[\varphi]=\varphi$. By \eqref{h}, \eqref{ave} and \eqref{eqgeo2}, we can rewrite $h_{\delta}$ as $h_{\delta}=h(\varphi_{\delta})$ whenever $\delta$ is small enough. Furthermore, we have 
\begin{equation}\label{vpes}
\lim\limits_{\delta\to0}\norm{\varphi_{\delta}}_{{\mathcal{W}}_0}=0.
\end{equation}  
\end{remark} 
\subsection{Proof of \autoref{main}}
Let us start this section with a basic fact. 
\begin{lemma}\label{anal}
The mapping defined by
\[
\Theta: \mathcal{F} \times \mathcal{W}_{\mathbb{R} }\to  \mathcal{W}_{\mathbb{R} },\ (g,\phi)\mapsto g\circ \phi
\]
is analytic.
\end{lemma}  
\begin{proof}
 For every $g\in\mathcal{F} $ and $n\geq0$, set 
 \[
 a_n(g) :=\Bigl(\frac{{\pi}^n}{n!}\Bigr)^{1/2}\inner{g}{e_n}.
 \]
 It is easy to see that $a_n$ is a linear functional on $\mathcal{F} $ whose norm is at most $(\frac{{\pi}^n}{n!})^{\frac{1}{2}} $ and 
 \[
 g(z)=\sum\limits_{n=0}^{\infty}a_n(g)z^n,\quad z\in\mathbb{C}.
 \]
 Thus, we have 
 \[
 \Theta(g,\phi)=\sum\limits_{n=0}^{\infty}a_n(g){\phi}^n,
 \]
 which yields the desired conclusion.
\end{proof}  
Let $F(\varphi_{\delta})=(f_1(\varphi_{\delta}),f_2(\varphi_{\delta}))\in\mathcal{F} ^2$ be the orthonormal system in \autoref{anaei} such that $\Omega=\Omega_{h_{\delta}}=\{u_F>t \}$ for some $t=t_{h_{\delta}}>0$. Here $f_j(\varphi) :=f_j(h(\varphi)), j=1,2$. Set 
\[
\Xi(\varphi)(\theta) :=u_{F(\varphi)}[(R+h(\varphi)(\theta))e^{i\theta}],\quad (\varphi,\theta)\in B_{\mathcal{W}_0}(0,\iota)\times[-\pi,\pi].
\]
Similarly, for $j=1,2$, we can also define 
\[
\Xi_j(\varphi)(\theta) :=u_{f_j(\varphi)}[(R+h(\varphi)(\theta))e^{i\theta}],\quad (\varphi,\theta)\in B_{\mathcal{W}_0}(0,\iota)\times[-\pi,\pi].
\]
One easily checks that $\Xi=\Xi_1+\delta \Xi_2$. By the definition of $\Omega_{h_{\delta}}$, we have $\Xi(\varphi_{\delta})(\theta)=t$ for all $\theta\in[-\pi,\pi]$, which implies that 
\[
\Pi(\varphi_{\delta}) :=P(\Xi(\varphi_{\delta}))=0.
\]
Here, we denote by $P$ the orthogonal projection from $\mathcal{W}_{\mathbb{R} }$ onto $\mathcal{W}_0$. In what follows, it suffices to prove that $\varphi=0$ is the unique zero of $\Pi$ in a small neiborhood of $0$ when $\delta$ is sufficiently small. For this reason, we need the following lemma:
\begin{lemma}\label{multi}
For $j=1,2$, set $\Pi_j :=P\circ\Xi_j$. Then $\Pi_j$ is an analytic mapping on $\mathcal{W}_0$. Furthermore, we have 
\begin{equation}\label{pi1}
\bigl((\mathcal{D}\Pi_1)(0)[\varphi]\bigr)^{\wedge}(n)=-2\pi Re^{-s}(1-\alpha_{|n|}(s))\widehat{\varphi}(n),\quad |n|\geq2.
\end{equation}  
Here, $\alpha_n(s)=\frac{s^n/n!}{\sum\limits_{j=1}^ns^j/j!}$.
\end{lemma}  
\begin{proof}
Set
\[
H(\varphi)(\theta) :=(R+h(\varphi)(\theta))e^{i\theta}.
\]
Then we obtain that 
\[
\Xi_j(\varphi)=\Theta(f_j(\varphi),H(\varphi))\overline{\Theta(f_j(\varphi),H(\varphi))}e^{-\pi|H(\varphi)|^2}.
\]
Thus, the analyticity follows from \autoref{anal}. In addition, since for every real-valued integrable function $f$, one has $\widehat{f}(n)=\overline{\widehat{f}(-n)}$, it suffices to prove \eqref{pi1} for $n\geq2$. By the definition of $f_1(h)$, we have 
\[
\mathcal{T}(h)f_1(h)=\lambda_1(h)f_1(h),
\]
which yields that 
\[
\Bigl(\mathcal{D}\mathcal{T}(0)[h]\Bigr)e_0+T_0\Bigl((\mathcal{D}f_1)(0)[h]\Bigr)=\Bigl((\mathcal{D}\lambda_1)(0)[h]\Bigr)e_0+\lambda_1(s)\Bigl((\mathcal{D}f_1)(0)[h]\Bigr)
\]
and hence 
\[
\inner{\Bigl(\mathcal{D}\mathcal{T}(0)[h]\Bigr)e_0}{e_n}+\inner{T_0\Bigl((\mathcal{D}f_1)(0)[h]\Bigr)}{e_n}=\lambda_1(s)\inner{(\mathcal{D}f_1)(0)[h]}{e_n}
\]
for every $n\geq1$. Since $\mathcal{T}(h)$ is a self-adjoint operator, we have
\[
\inner{T_0\Bigl((\mathcal{D}f_1)(0)[h]\Bigr)}{e_n}=\inner{(\mathcal{D}f_1)(0)[h]}{T_0e_n}=\lambda_{n+1}(s)\inner{(\mathcal{D}f_1)(0)[h]}{e_n}.
\]
Therefore, we obtain that 
\[
\inner{(\mathcal{D}f_1)(0)[h]}{e_n}=\frac{1}{\lambda_1(s)-\lambda_{n+1}(s)}\inner{\Bigl(\mathcal{D}\mathcal{T}(0)[h]\Bigr)e_0}{e_n},\quad n\geq1.
\]
Combining the identity above with \eqref{eqop} and \autoref{inner}, we have 
\begin{equation}\label{df1}
(\mathcal{D}f_1)(0)[h]=\sum\limits_{n=1}^{\infty}\frac{2\pi Re^{-s}}{\lambda_1(s)-\lambda_{n+1}(s)}\cdot\Bigl(\frac{s^n}{n!}\Bigr)^{\frac{1}{2}}\widehat{h}(n)e_n.
\end{equation}
Thus, we directly derive \eqref{pi1} from \eqref{hpro}, \eqref{df1}, and the chain rule.
\end{proof}  
\begin{proof}[Proof of \autoref{main}]
It is easy to check that $\inf\limits_{n\geq2}(1-\alpha_n(s))>0$. By \eqref{pi1}, we conclude that $(\mathcal{D}\Pi_1)(0)$ is an invertible operator on $\mathcal{W}_0$. Thus, there exists $\kappa>0$ such that every operator $T$ on $\mathcal{W}_0$ satisfying 
\begin{equation}\label{es1}
\norm{T-(\mathcal{D}\Pi_1)(0)}_{\mathcal{W}_0\to\mathcal{W}_0}<\kappa
\end{equation}
is invertible. Set \[\delta_1 :=\frac{\kappa}{10(\norm{(\mathcal{D}\Pi_2)(0)}_{\mathcal{W}_0\to\mathcal{W}_0}+1)}.\] Then by \eqref{es1}, for all $\delta\in[0,\delta_1)$ we obtain that $(\mathcal{D}\Pi)(0)$ is an invertible operator on $\mathcal{W}_0$. Therefore, we use the vector-valued version of the inverse function theorem \cite[Theorem 15.2]{MR787404} to conclude that there exists $\tau>0$ such that $\varphi=0$ is the unique zero of $\Pi$ in $B_{\mathcal{W}_0}(0,\tau)$. However, by \eqref{vpes}, we can choose $\delta_2>0$ such that $\delta\in[0,\delta_2)$ implies $\varphi_{\delta}\in B_{\mathcal{W}_0}(0,\tau)$. Setting $\delta_s :=\min\{\delta_1,\delta_2\}$, we conclude the proof.
\end{proof}  

\section*{Acknowledgments}
The author is grateful to his advisor João P. G. Ramos for introducing him to this
project and for several inspiring comments. The author also thanks Jorge Isaac Rubiano for many helpful discussions.

\section*{Statement of LLM usage}
Large language models played a constructive role in the preparation of this article. As for the mathematical content, the main ideas of this paper are of human origin. GPT 6.1 was used to deal with some technical issues. In addition, the discussions with it provided constructive ideas and valuable inspiration. Meanwhile, the author simplified and revised many of the proofs. In other respects, the author used GPT 6.1 to resolve some issues encountered while editing LaTeX, and used Gemini and Doubao to improve the English wording. The author takes full responsibility for all the content of this paper.
 \printbibliography
\par\medskip
\begingroup
  \normalfont\fontsize{10pt}{12pt}\selectfont
  \setlength{\parindent}{1em}
  \indent (Yifan Zhang)
  \textsc{Instituto de Matemática Pura e Aplicada, Rio de Janeiro, Brazil}\par
  \indent\textit{Email address:}\enspace
  {\urlstyle{tt}\nolinkurl{yifan.zhang@impa.br}}\par
\endgroup
\end{document}